\RequirePackage{amsmath}
\documentclass[runningheads]{llncs}
\usepackage{microtype}
\usepackage[T1]{fontenc}
\usepackage{graphicx,enumitem}
\usepackage{subfig}
\usepackage{booktabs}

\usepackage[utf8]{inputenc}

\usepackage{amssymb}
\usepackage{xspace}
\usepackage{xcolor}
\usepackage[
backend=biber,
style=ieee,
sorting=ynt
]{biblatex}
\usepackage{authblk}

\usepackage{color}
\title{Conditional traffic-like rules for particle-flow simulation in cellular automata}

\author{
Barbara Wolnik\inst{1,3} \and 
Dominik Michał Falkiewicz\inst{5} \and 
Witold Bołt\inst{2} \and 
Adam Rutkowski\inst{4} \and 
Bernard De Baets\inst{3}
}

\institute{
Institute of Mathematics, Faculty of Mathematics, Physics and Informatics, University of Gdańsk, Wita Stwosza 57, 80-308 Gdańsk, Poland\\ 
\email{barbara.wolnik@ug.edu.pl}
\and
Jit Team, Łużycka 8C, 81-537 Gdynia, Poland\\ 
\email{witold.bolt@jit.team}
\and
KERMIT, Department of Data Analysis and Mathematical Modelling, Faculty of Bioscience Engineering, Ghent University, Coupure links 653, B-9000 Gent, Belgium\\ 
\email{bernard.debaets@ugent.be}
\and 
Institute of Theoretical Physics and Astrophysics, Faculty of Mathematics, Physics and Informatics, University of Gdańsk, Wita Stwosza 57, 80-308 Gdańsk, Poland\\ 
\email{adam.rutkowski@ug.edu.pl}
\and  
Faculty of Mathematics, Physics and Informatics, University of Gdańsk, Wita Stwosza 57, 80-308 Gdańsk, Poland\\ 
\email{dominik.falkiewicz@phdstud.ug.edu.pl}
}

\date{\today}

\authorrunning{B. Wolnik, D. M. Falkiewicz, W. Bołt, A. Rutkowski, B. De Baets}

\titlerunning{Conditional traffic-like rules for particle-flow simulation in cellular automata}

\newcommand{\ff}[6]{
{f}\!\left(\!
\begin{array}{ccc}
#1 & #2 & #3 \\
#4 & #5 & #6
\end{array}\!
\right)
}
\newcommand{\Id}{\mathrm{Id}}
\newcommand{\Traffic}[1]{\mathrm{T#1}}
\newcommand{\Shift}[1]{\mathrm{S#1}}
\newcommand{\Cset}[1]{C_{#1}}
\newcommand{\C}{C}

\newcommand{\bx}{\mathbf{x}}

\newcommand{\cbi}{\mathbf{i}}
\newcommand{\bj}{\mathbf{j}}

\newcommand{\Z}{\mathbb{Z}}

\newcommand{\di}{\operatorname{dist}}

\newcommand{\vv}[2]{{\overrightarrow{{\bf #1}_{#2}}}}
\newcommand{\smallnwarrow}{\scalebox{0.5}{$\nwarrow$}}
\newcommand{\smallnearrow}{\scalebox{0.5}{$\nearrow$}}
\newcommand{\smallswarrow}{\scalebox{0.5}{$\swarrow$}}

\newcommand{\eg}{\textit{e.g.}}

\newcommand{\ie}{\textit{i.e.}}

\usepackage{letltxmacro}
\LetLtxMacro{\originaleqref}{\eqref}
\renewcommand{\eqref}{Eq.\ \originaleqref}

\newcommand{\DD}{\textit{D}}

\begin{document}

\maketitle

\begin{abstract}
This paper presents a novel approach to the description and understanding of two-dimensional binary cellular automata with the Moore neighborhood that preserve the number of active cells. Such dynamical systems are known to successfully capture particle dynamics and are often used as theoretical models of physical processes where some conservation laws have to be taken into account. Unfortunately, to date, there are no tools for designing such non-trivial cellular automata or for studying their properties, not to mention finding them all and describing their dynamics (even the order of magnitude of their number is unknown). We believe that the novel framework unfolded in this paper will make it possible to overcome all these challenges.

\end{abstract}

\section{Introduction}
Due to their ability to represent complex systems using simple local rules, cellular automata (CAs) are widely used to model various natural phenomena (see \cite{PhysRevE.111.025404, PhysRevE.111.054306, PhysRevE.111.054118, PhysRevE.111.014419} for recent applications). In particular, they have been used for many decades for modeling all kinds of physical systems, especially those where different conservation laws have to be taken into account (\cite{Chopard2009, Goles1992, Vichniac1984}). One example of such a system is a set of moving particles. In typical scenarios, it is usually assumed that the particles under consideration are indistinguishable and move in a hypothetical, $d$-dimensional regular grid, while their number is constant (\ie, the particles cannot disappear or appear out of nowhere). It is also often assumed that there can be at most one particle in each cell at any point in time (reducing the size of the cells if necessary).
In CA modeling, the combination of these assumptions leads to so-called number-conserving binary CAs.

When designing a number-conserving binary CA that models a given phenomenon, one starts by extracting the very local behavior of the particles, which is then translated into the language of the local rule of the created CA. The simplest and the most well-known example of this procedure are the so-called \emph{traffic rules}. Here, we start with a one-dimensional grid of cells in which each particle moves in a fixed direction (say, to the right) but only if the cell in front of it is empty. The described behavior is sufficient to construct a CA with a very simple local rule of radius one, which can be represented in the following tabular form called LUT (LookUp Table):
\begin{table}[ht]
	\centering
		\begin{tabular}{|l|c|c|c|c|c|c|c|c|}
		\hline
		 \rule{0pt}{1.2em} neighborhood pattern  & 111 & 110 & 101 & 100 & 011 & 010 & 001 & 000   \\
		\hline
		 \rule{0pt}{1.2em} output of the local rule &	$1$ & $0$ & $1$ & $1$ & $1$ & $0$ & $0$ & $0$ \\
		 \hline
	    \end{tabular}
	\label{tab:ECA}
\end{table}

Nowadays, when so many CAs can be generated by computers in terms of LUTs (for example, using advanced genetic algorithms), a quite natural question arises: is the converse procedure possible? That is, given a number-conserving binary CA, are we able to describe it using a particle displacement representation?
It turns out that this task, in general, is not easy.

First of all, it should be noted that for a given binary CA in terms of its LUT, it is possible to check whether it is number-conserving or not. In mathematical terms, this boils down to checking a certain system of linear equations (see, for instance, \cite{Durand2003}), so it is not too complicated and easy to outsource to a computer.
The last point is key because the number of columns in LUT for a binary local rule with a neighborhood of size $m$ equals $2^m$. So, for example, considering a two-dimensional binary CA with the Moore neighborhood (consisting of 9 cells), we get a LUT with as many as 512 values. 
With only these 512 LUT values, it is difficult to see, and hence also to understand, what the dynamics of a given CA looks like. Most eagerly, we would like to use some algorithm and, with the help of a computer, answer the question of whether or not a given CA has some property we are interested in. 
Unfortunately, most properties of interest inspired by physical phenomena are undecidable (often in dimensions greater than one \cite{kari2005, kutrib2008cellular}). This is the case, for example, with reversibility: in general, there are no algorithms to determine if a given two-dimensional CA is injective or surjective~\cite{kari1996representation}.
However, it is believed that a significant portion of properties that are undecidable for general CAs may be decidable for the subclass consisting of number-conserving CAs.
Thus, translating the LUT of a given number-conserving binary CA in terms of particle movement gives us an opportunity to understand how this CA works ``at the microscopic level'' and analyze it from a new perspective.

As for one-dimensional number-conserving CAs, it can be said that the task of translating a given CA into a particle displacement representation has been quite well studied (and not only for the binary set of states). 
It was proven that every one-dimensional number-conserving CA can be interpreted using motion representations; moreover, any such CA can be uniquely characterized by some kind of \emph{canonical form} of motion representation (for details, see \cite{BF98,BoccaraF02,Pivato,fuks-2000-class,MOREIRA2004285}).
Another very interesting characterization of one-dimensional number-conserving CAs was realized through the use of \emph{flow functions}, which describe how many particles cross the boundary between two adjacent cells, depending on the neighborhood of this boundary~\cite{Redeker22}. 

In the two-dimensional case, the topic of motion representation has not been well explored, although in this case it was also shown that the dynamics of any two-dimensional number-conserving CA can be expressed in terms of particle displacements~\cite{KariTaati}. However, for such a CA usually there is an infinite number of motion representations, and it is hard to define a ``canonical'' one. 

When it comes to higher dimensions, there are practically no results on motion representation of number-conserving CAs, even in the binary case. Although it is rather clear that the results from~\cite{KariTaati} can be generalized to three or more dimensions, research in this direction has not been continued at all. 
The exception to this are number-conserving CAs with the von Neumann neighborhood.
First, there is a characterization of such CAs in terms of some flow functions in the vertical, horizontal, and diagonal direction~\cite{TI} (but again the same problem occurs: the flow functions are not uniquely defined).
However, if we focus solely on the binary case, then using another mathematical tool, the split-and-perturb decomposition~\cite{decomposition}, constructed specifically for the study of $d$-dimensional number-conserving CAs with the von Neumann neighborhood, we can easily find all binary such CAs.
And then it turns out that there are only $4d+1$ of them: the identity rule, the shift rules and the traffic rules in each of the $2d$ directions~\cite{PhysRevE.100.022126}.
In this way, the structure of the $d$-dimensional number-conserving binary CAs with the von Neumann neighborhood is fully unveiled: regardless of the dimension~$d$, all of these CAs are trivial, as they are intrinsically one-dimensional, which makes it possible to describe the motion representation of such automata in detail~\cite{PhysRevE.100.022126}. This result is due to the geometric properties of the von Neumann neighborhood (see~\cite{decomposition}, for more details).

Unfortunately, this approach is not possible in the case of the Moore neighborhood. Already for $d=2$, the cardinality of the set of all binary CAs with such a neighborhood reaches the astronomical value of $2^{512}$ and no one has yet managed to answer the question of how many among them are number-conserving and what their dynamics look like (using the aforementioned necessary and sufficient conditions~\cite{Durand2003} is beyond the computing power of computers).
The only results known to the authors regarding approaches to finding and describing all two-dimensional number-conserving binary CAs with the Moore neighborhood can be found in~\cite{7424749}, where an algorithm is proposed that allows one to generate all two-dimensional ``motion-representable'' number-conserving binary CAs defined on $n\times m$ neighborhoods. It turns out that in the cases $2\times 2$ and $2\times 3$, this algorithm effectively searches for all number-conserving binary CAs for the neighborhood under consideration. However, due to its complexity, it has never been used in case of the Moore neighborhood for $d=2$, not to mention larger values.

In this paper, we propose a very new way of looking at number-conserving binary CAs with the Moore neighborhood as ``modular laws'' governing the movement of particles. In this new insight, each such law consists of two parts: an overarching regulation of the movement of the entire system in a specified direction $\Omega$ and a set $\Lambda$ of regulations that are applied locally. 
The representation of the CA as a couple $(\Omega,\Lambda)$ has the following advantages. First of all, both parts are uniquely determined. Secondly, instead of a long mathematical description in the form of a LUT, we obtain a short, transparent explanation in the form of particle movement. 
Moreover, this representation allows us, as we will see later in the paper, to track the movement of each particle separately.
This makes it easier, among other things, to understand the dynamics of the CA under study and to verify its properties, like, for example, reversibility.

It turns out that our approach makes it easy to describe all two-dimensional number-conserving binary CAs both with the von Neumann neighborhood and with the $2\times 3$-neighborhood. At this point, it is not known whether it also allows us to describe all two-dimensional number-conserving binary CAs with the Moore neighborhood. However, it is already a very useful tool for designing highly non-trivial such CAs and studying their properties.

The remainder of this paper is organized as follows. The key definitions and a description of the new approach are presented in Section~2. In Section~3, the proposed approach is used to find all two-dimensional number-conserving binary CAs with the $2\times 3$-neighborhood and describe their dynamics. Section~4 contains an example application of the introduced approach in the case of the Moore neighborhood. Conclusions and future work are reported in Section 5.

\section{Cellular automata with radius one that simulate particle flow}

In this paper, we assume that the cellular space considered $\C$ is infinite and corresponds to the whole space $\Z^d$, where $d$ is a positive integer. However, the results are valid also for any finite cellular space that corresponds to
some $d$-dimensional cuboid of size ${\bf N}=(n_1,n_2,\ldots,n_d)$ with periodic boundary conditions, \ie,
\arraycolsep=2pt
\begin{eqnarray*}
\C &=& \C_{\bf N}=\left(\Z / n_1\Z\right)\times \left(\Z / n_2\Z\right)\times\cdots\times \left(\Z / n_d\Z\right)\\
   &=& \{0,1,\ldots,n_1-1\}\times\{0,1,\ldots,n_2-1\}\times\cdots \times\{0,1,\ldots,n_d-1\} 
\end{eqnarray*}
(provided that $n_1,n_2,\ldots,n_d> 2$).
We use $Q=\{0,1\}$ as the set of states, where the state $1$ is interpreted as a particle and the state $0$ as an empty space.

\begin{figure}[!ht]
\centering
 \subfloat[]{
\includegraphics[height=0.30\textwidth]{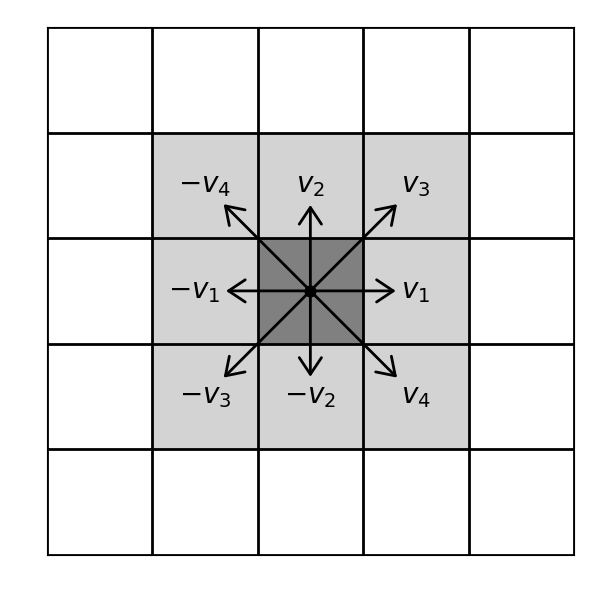}
}
 \subfloat[]{
\includegraphics[height=0.30\textwidth]{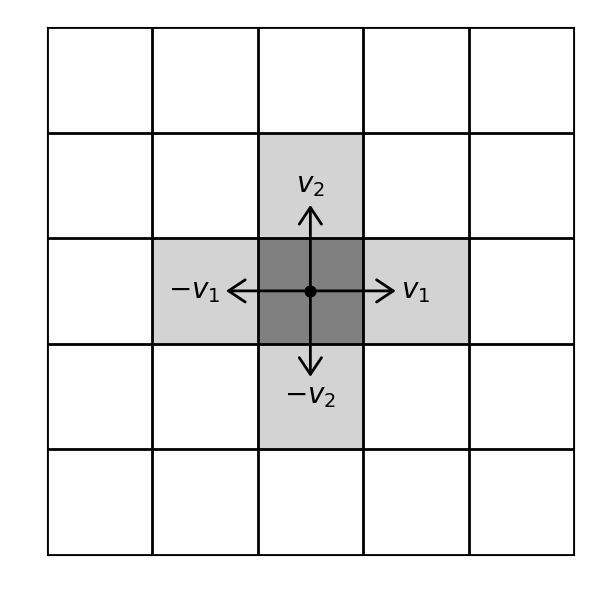}
}
 \subfloat[]{
\includegraphics[height=0.30\textwidth]{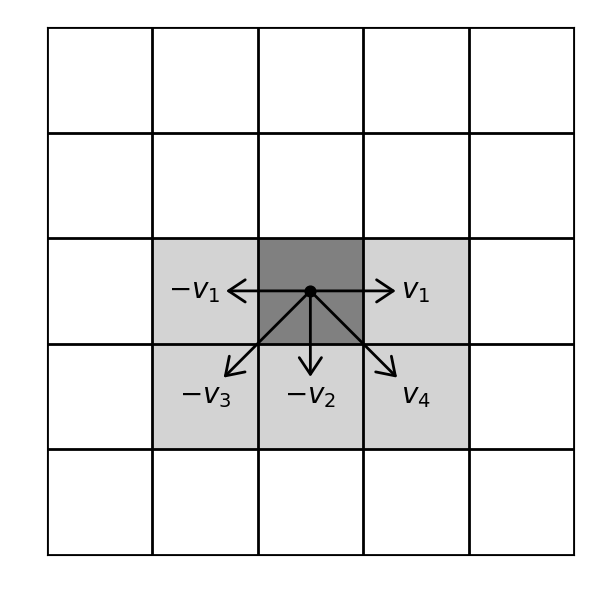}
}
\caption{The considered two-dimensional neighborhoods of cell $\cbi$ (a)  the  Moore neighborhood, (b) the von Neumann neighborhood, (c) the $2\times3$-neighborhood. The cell $\cbi$ is marked by dark gray and the cells in its neighborhood are marked by light gray.}
\label{neighborhoodFig}
\end{figure}

The neighborhood of a cell is usually given by a vector of vectors $V=(\vv{v}{1},\vv{v}{2},\ldots, \vv{v}{m})$, which means that, for a cell $\cbi\in\C$, the cells $\cbi+\vv{v}{1},\cbi+\vv{v}{2},\ldots,\cbi+ \vv{v}{m}$ are its neighbors (it is usually assumed that the zero vector belongs to $V$, \ie, that any cell is its own neighbor). Among the many different ones that could be adopted (see, \eg,~\cite{Zaitsev}), the two most popular ones are the Moore neighborhood and the von Neumann neighborhood (see Figure~\ref{neighborhoodFig}(a)(b)). Both can be seen as closed unit balls centered on the cell in question, but they differ in terms of the metric considered in $\Z^d$: the Moore neighborhood is defined by the Chebyshev distance
\begin{equation}
\di_\mathrm{Ch} (\cbi ,\bj )=\max_{1\leq k\leq d} |i_k-j_k|\,,
\end{equation}
whereas the von Neumann neighborhood is defined by the Manhattan distance:
\begin{equation}\label{di}
\di_\mathrm{M}(\cbi ,\bj )=\sum_{k=1}^{d} |i_k-j_k|\,,
\end{equation}
where $\cbi =(i_1,i_2,\ldots,i_d)\in\C$ and $\bj =(j_1,j_2,\ldots,j_d)\in\C$. Thus, in the case of two dimensions, the Moore neighborhood is given by \[
V_M=(-\vv{v}{4},\vv{v}{2},\vv{v}{3},-\vv{v}{1},\vv{v}{0},\vv{v}{1},-\vv{v}{3},-\vv{v}{2}, \vv{v}{4})\,,
\]
whereas the von Neumann neighborhood is given by
\[
V_N=(\vv{v}{2},-\vv{v}{1},\vv{v}{0},\vv{v}{1},-\vv{v}{2}),
\]
where $\vv{v}{0}=[0,0]$, $\vv{v}{1}=[1,0]$, $\vv{v}{2}=[0,1]$, $\vv{v}{3}=[1,1]$, and $\vv{v}{4}=[1,-1]$.
However, if we want to consider a two-dimensional CA with radius one, we can use any non-empty subset of $V_M$ as the neighborhood, and for our purposes we will also consider the asymmetric six-cell neighborhood given by 
\[
V=(-\vv{v}{1},\vv{v}{0},\vv{v}{1},-\vv{v}{3},-\vv{v}{2}, \vv{v}{4})\,,
\]
which we will call the $2\times 3$-neighborhood (see Figure~\ref{neighborhoodFig}(c)).

A \emph{configuration} is any mapping from the grid $\C$ to $Q$ and the set of all configurations is denoted by $X$. 
The state of a cell $\cbi\in\C$ in a configuration $\bx\in X$ is denoted by $x_\cbi$. A configuration $\bx\in X$ is called \emph{finite} if at most finitely many cells are assigned a non-zero state, \ie, the set $\{\cbi\in \C\mid x_\cbi\neq 0\}$ is finite. 

A binary CA acts according to some \emph{local rule} $f:\{0,1\}^m\to \{0,1\}$: the state of the cell $\cbi$ in the next time step is given by $f(x_1,x_2,\ldots,x_m)$, where $x_1,x_2,\ldots,x_m$ are the current states of the cells $\cbi+\vv{v}{1},\cbi+\vv{v}{2},\ldots,\cbi+ \vv{v}{m}$, respectively. 
A given local rule $f$ induces a \emph{global rule} $F: X\to X$ defined for any $\bx\in X$ and $\cbi\in\C$ as follows
\[
F(\bx)_\cbi=f(x_{\cbi+\vv{v}{1}},x_{\cbi+\vv{v}{2}},\ldots,x_{\cbi+ \vv{v}{m}})\,.
\]

We are interested in binary CAs that can be regarded as a model of a $d$-dimensional system of indistinguishable particles, in which the particles can neither disappear nor appear, but only move within the grid. These assumptions led to the study of \emph{number-conserving} binary CAs, for which the sum of the states is the same at any time step (or, equivalently, the number of 1s does not change during the evolution of the system).

\begin{definition}
A CA with global rule $F$ is number-conserving if for every finite configuration $\bx \in X$ it holds that $\sum_{\cbi \in\C}F(\bx)_\cbi = \sum_{\cbi \in\C}\bx_\cbi$.  
\end{definition}

In the one-dimensional case, the Moore neighborhood and the von Neumann neighborhood coincide and lead to $256$ Elementary CAs (ECAs). Among them, there are only five number-conserving ones: the identity rule ($\Id$), the shift-left rule ($\Shift{l}$), the shift-right rule ($\Shift{r}$), the traffic-left rule ($\Traffic{l}$), and the traffic-right rule ($\Traffic{r}$). Most often, these rules are interpreted as follows. The identity rule and the shift rules are stationary: the pattern of a configuration does not change, it can only stand motionless or move to the right or to the left. 
More interesting are the traffic rules, which are often used as the simplest models of road traffic flow. The states $1$ and $0$ are then interpreted as `car' and `empty space', respectively. In such a model, cars move only if there is an empty space in front of them (see, e.g.,~\cite{PhysRevE.55.R2081}).
As we can see, the identity rule and the shift rules are basically the same in this interpretation.

In this paper, we provide a slightly different perspective on number-conserving binary CAs and look at them as `laws' governing the movement of particles. As mentioned in the introduction, in this new insight, each such law consists of two parts: an overarching regulation of the movement of the entire system in a specified direction (the identity or a shift in some fixed direction), denoted by $\Omega$, and a set of regulations that are applied locally (conditional traffics), denoted by $\Lambda$.
In this new approach, the one-dimensional number-conserving binary CAs described above have the following new interpretation in terms of $(\Omega, \Lambda)$: the shift-left rule corresponds to $(\Shift{l}, \emptyset)$ and the shift-right rule corresponds to $(\Shift{r}, \emptyset)$. Moreover, the shift rules do not allow any locally applied regulations. The situation is different when the overarching regulation is the identity rule, because then the traffic rules are allowed as locally applied regulations, and finally the identity rule, the traffic-left rule and the traffic-right rule correspond to $(\Id, \emptyset)$, $(\Id, \{\Traffic{l}\})$ and $(\Id, \{\Traffic{r}\})$, respectively.

As we mentioned earlier, in the case of the von Neumann neighborhood, it has been shown that increasing the dimension does not, unfortunately, result in the appearance of other types of number-conserving binary CAs. Indeed, there are exactly $4d+1$ $d$-dimensional number-conserving binary CAs with the von Neumann neighborhood: the identity rule and the shift and traffic rules in each of the $2d$ possible directions (see~\cite{PhysRevE.100.022126}).   
Using our new approach, we can describe the set of all $d$-dimensional number-conserving binary CAs with the von Neumann neighborhood as follows: 
\begin{itemize}
    \item Class $(\Omega,\emptyset)$ consisting of $2d+1$ number-conserving binary CAs with $\Lambda = \emptyset$: the identity rule and the shift rule in each of $2d$ possible directions,
    \item Class $(\Id,\Lambda)$ consisting of $2d$ number-conserving binary CAs with $\Omega=\Id$ and $\Lambda$ being the traffic rule in each of $2d$ possible directions.
\end{itemize}

In the next section, we will show how our new approach to understanding number-conserving binary CAs allows us to find all two-dimensional such CAs with the $2\times 3$-neighborhood and provide a unique description in the form of couples $(\Omega,\Lambda)$.

\section{Two-dimensional number-conserving binary CAs with the $2\times 3$-neighborhood}

In this section we introduce a generalization of traffic rules and describe the possibility of ``mixing'' such rules. Since, as we will show, the $2\times 3$-neighborhood allows traffic to occur even when an overarching regulation is one of the shift rules (namely, the shift-up rule $\Shift{u}$), we consider two options for $\Omega$: the identity rule and the shift-up rule. Although it is not necessary, we consider these cases separately for the convenience of the readers. The approach based on $\Omega$ and $\Lambda$ allows us to find a set of two-dimensional number-conserving binary CAs with $2\times 3$-neighborhood.
Then, we show that all two-dimensional number-conserving binary CAs with this neighborhood can be generated in this way.

For convenience, we will use the following simplifications. For two cells $\cbi,\cbi'\in\C$ we will write that $\cbi$ sees $\cbi'$ if cell $\cbi'$ belongs to the neighborhood of cell $\cbi$. Moreover, we will use $(i,j)$ to refer to a given cell instead of $\cbi$.

\subsection{Conditional traffic rules in the case $\Omega=\Id$}
In the one-dimensional case the traffic-right rule acts as follows: If there is a particle in cell $i$, while cell $i+1$ is empty, then in the next time step the particle from cell $i$ moves into cell $i+1$. Such a regulation is possible because cell $i+1$ belongs to the neighborhood of cell $i$ and vice versa. This gives both cells the information they need: cell $i$ sees the empty space on the right, so it knows it will get rid of the particle, while cell $i+1$ sees the particle on the left, so it knows it will get it.
In the case of the $2\times 3$-neighborhood, the extra dimension allows that such kind of regulation may be more complicated.

\begin{figure}[!ht]
\centering
\includegraphics[height=0.3\textwidth]{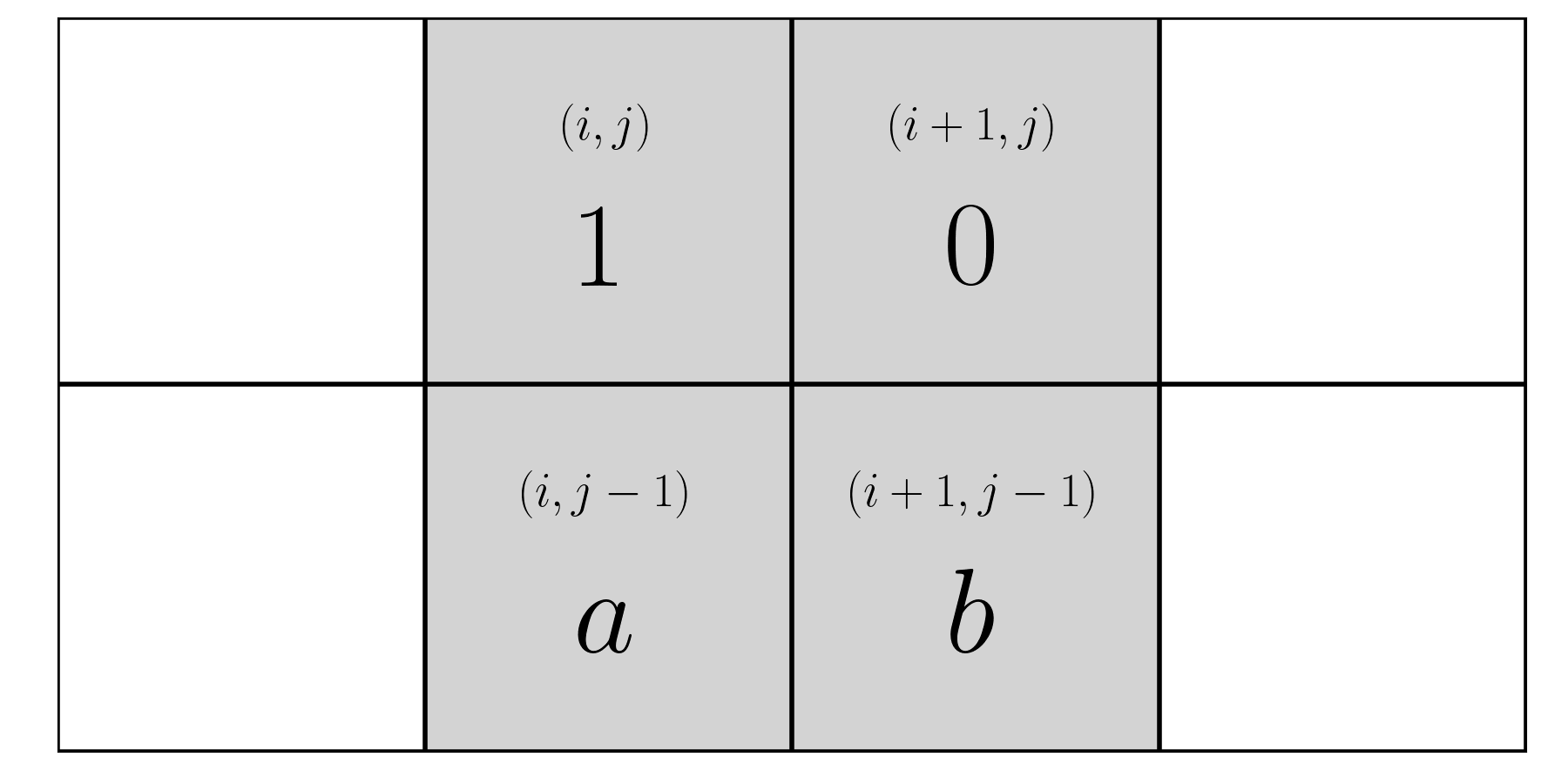}
\caption{Under the traffic-right rule, the particle from cell $(i,j)$ will always move into cell $(i+1,j)$ as long as $(i+1,j)$ is empty. Conditional traffic rules can be defined, where the movement only happens depending on the values of $a$ and $b$.}
\label{fig:TR}
\end{figure}

In fact, there are up to four cells that are seen both by $(i,j)$ and its right neighbor $(i+1,j)$: $(i,j)$, $(i+1,j)$, $(i,j-1)$, and $(i+1,j-1)$. So, if there is a particle in cell $(i,j)$ and cell $(i+1,j)$ is empty, whether the particle moves to the right may depend on what is in cells $(i,j-1)$ and $(i+1,j-1)$ (see Figure~\ref{fig:TR}). 
Such type of regulation governing the movement of particles will be called a \emph{conditional traffic-right rule}.

\begin{definition}
Let $C$ be a non-empty subset of $\{00,01,10,11\}$. The conditional traffic-right rule under condition~$C$ acts as follows. If $\bx\in X$ satisfies $x_{i,j}=1$ and $x_{i+1,j}=0$, then $F(\bx)_{i,j}=0$ and $F(\bx)_{i+1,j}=1$ if and only if $x_{i,j-1}x_{i+1,j-1}\in C$.
\end{definition}
In other words, if there is a particle in cell $(i,j)$ and cell $(i+1,j)$ is empty, then the particle moves to the right if and only if the pattern of cells $(i,j-1)$ and $(i+1,j-1)$ belongs to $C$.
The conditional traffic-right rule under condition $C$ is denoted by $\Traffic{r}_{C}$.
Since there are exactly $15$ non-empty subsets of $\{00,01,10,11\}$, we have $15$ conditional traffic-right rules. Note that if $C=\{00,01,10,11\}$, then $\Traffic{r}_{C}$ is the usual traffic-right rule, \ie, the particle moves to the right whenever it sees an empty space in front of it.

We define conditional traffic-left rules ($\Traffic{l}_{C}$) in an analogous way.

\begin{definition}
Let $C$ be a non-empty subset of $\{00,01,10,11\}$. The conditional traffic-left rule under condition $C$ acts as follows. If $\bx\in X$ satisfies $x_{i,j}=1$ and $x_{i-1,j}=0$, then $F(\bx)_{i,j}=0$ and $F(\bx)_{i-1,j}=1$ if and only if $x_{i-1,j-1}x_{i,j-1}\in C$.
\end{definition}

This means that if there is a particle in cell $(i,j)$ and cell $(i-1,j)$ is empty, then the particle moves to the left if and only if the pattern of cells $(i-1,j-1)$ and $(i,j-1)$ belongs to $C$.  
Obviously, there are also exactly $15$ conditional traffic-left rules.

\begin{figure}[!ht]
\centering
 \subfloat[]{
\includegraphics[height=0.18\textwidth]{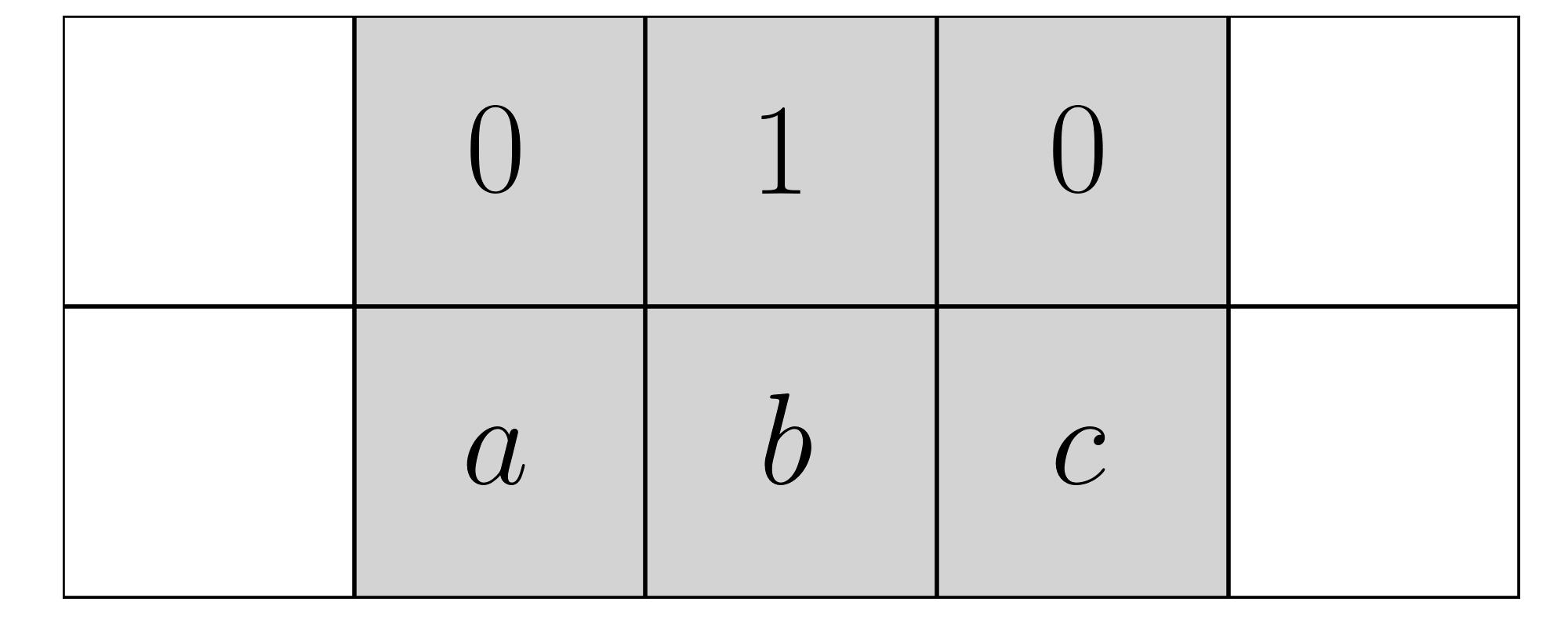}
}
\subfloat[]{
\includegraphics[height=0.18\textwidth]{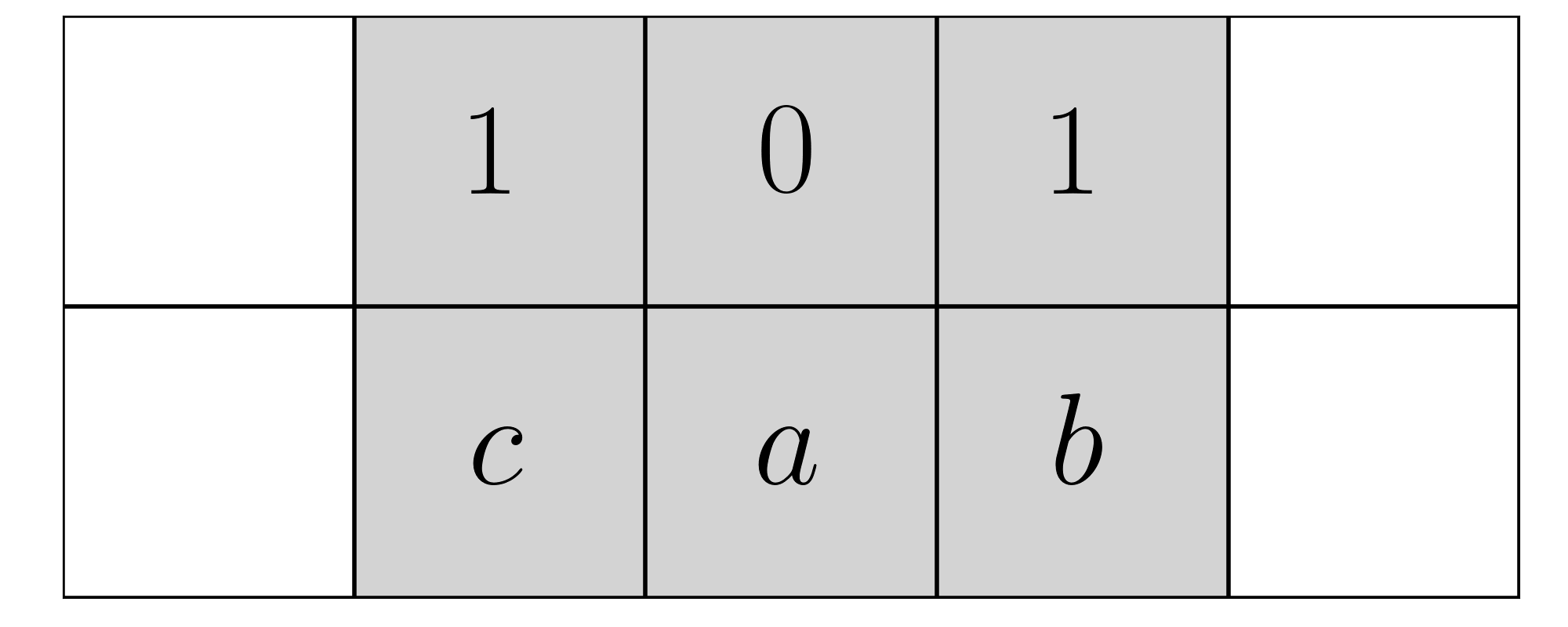}
}
\caption{Illustration of the requirements (R1) and (R2) for two conditional traffic rules $\Traffic{l}_{\Cset{l}}$ and $\Traffic{r}_{\Cset{r}}$ to co-exist within the same number-conserving binary CA. (a) If the pattern $ab$ allows 1 go to the left ($ab\in \Cset{l}$), then the pattern $bc$ cannot allow 1 go to the right ($bc\not\in \Cset{r}$).  (b) If the pattern $ab$ allows 1 go to the left ($ab\in \Cset{l}$), then the pattern $ca$ cannot allow 1 go to the right ($ca\not\in \Cset{r}$).}
\label{fig:R1R2}
\end{figure}

Note that the $2\times 3$-neighborhood does not allow any traffic-down rule at all, because while cell $(i,j)$ sees cell $(i-1,j)$, cell $(i-1,j)$ does not see cell $(i,j)$. For the same reason, traffic-up rules or traffic rules in both diagonal directions are not possible.
However, it is possible for two conditional traffic rules $\Traffic{l}_{\Cset{l}}$ and $\Traffic{r}_{\Cset{r}}$ to co-exist within the same number-conserving binary CA. Indeed, it is possible to choose conditions $\Cset{l}$ and $\Cset{r}$ in such a way that the following two requirements are met:
\begin{itemize}[leftmargin=1cm]
    \item[(R1)] No particle can move in two different directions at the same time.
    \item[(R2)] No empty cell can receive a particle from two different directions at the same time.
\end{itemize}
Both the situations described in (R1) and (R2) are shown in Figure~\ref{fig:R1R2}.

To ensure that (R1) is achieved, if $ab\in \Cset{l}$, then $bc\not\in \Cset{r}$, as shown in Figure~\ref{fig:R1R2}(a). And to ensure that (R2) is satisfied, if $ab\in \Cset{l}$, then $ca\not\in \Cset{r}$, as shown in Figure~\ref{fig:R1R2}(b).
In conclusion, if $ab\in \Cset{l}$, then $\bullet a, b\bullet\not\in \Cset{r}$ (the symbol $\bullet$ refers to any value). Therefore, there are only four possibilities for $\Cset{l}$ and $\Cset{r}$:
\begin{equation}\label{cases}
        \begin{array}{l}
        \Cset{l}=\{00\} \quad {\rm and }\quad \Cset{r}=\{11\}  \\
        \Cset{l}=\{01\} \quad {\rm and }\quad \Cset{r}=\{01\} \\
        \Cset{l}=\{10\} \quad {\rm and }\quad \Cset{r}=\{10\} \\
        \Cset{l}=\{11\} \quad {\rm and }\quad \Cset{r}=\{00\} \,.
    \end{array}
\end{equation}

In summary, if the overarching regulation is the identity, then, except for $\emptyset$, there are $34$ sets of locally applied regulations based on traffic rules: $15$ of the form $\{\Traffic{l}_{C}\}$, $15$ of the form $\{\Traffic{r}_{C}\}$, and four of the form $\{\Traffic{l}_{\Cset{l}}, \Traffic{r}_{\Cset{r}}\}$.

\subsection{Conditional traffic rules in the case $\Omega=\Shift{u}$}
Since any cell has five neighbors (excluding itself) in the case of the $2\times 3$-neighborhood, as many as five shift rules are possible: the shift-left rule, the shift-right rule, the shift-up rule (along the direction $\vv{v}{2}=[0,1]$), and two shift rules along the diagonal directions (\ie, along the vectors $\vv{v}{3}=[1,1]$ and $-\vv{v}{4}=[-1,1]$).

As already mentioned, the one-dimensional case does not allow for any non-empty set of locally applied regulations when the overarching regulation is any of the shift rules.
The same is true for the $2\times 3$-neighborhood when considering the shift-right, the shift-left or the shift rule along any of the diagonal directions.
However, if the overarching regulation is the shift-up rule (denoted by $\Shift{u}$), then the set of locally applied regulations can be non-empty.

\begin{figure}[!ht]
\centering
\includegraphics[height=0.3\textwidth]{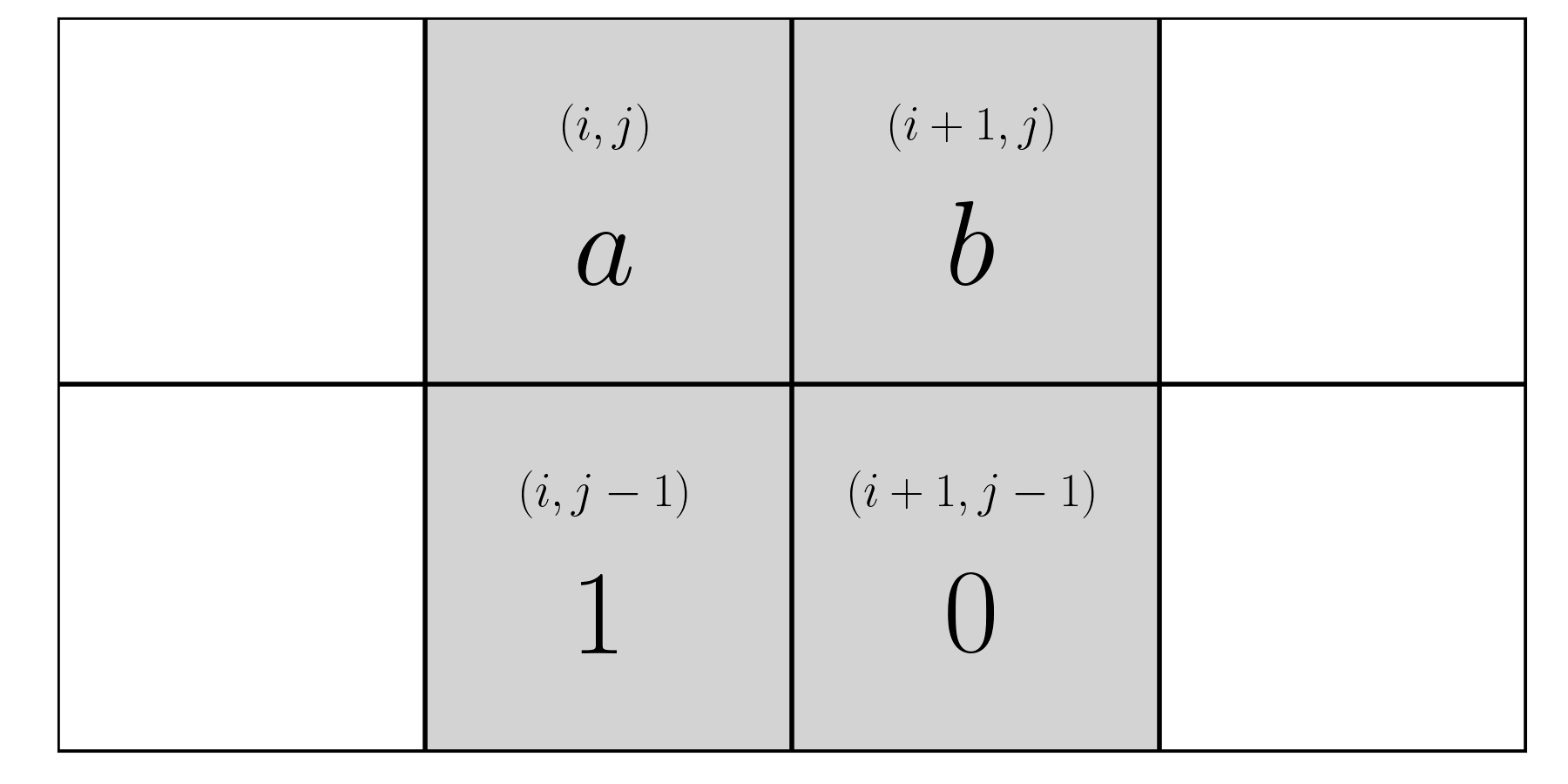}
\caption{Illustration of a shifted conditional traffic-right rule. Certain patterns $ab$ can allow the particle from cell $(i,j-1)$ to move to cell $(i+1,j)$ instead of $(i,j)$.}
\label{fig:shift-up}
\end{figure}

Indeed, cells $(i,j)$ and $(i+1,j)$, in the next time step,  receive the particles that are now in cells $(i,j-1)$, and $(i+1,j-1)$ (if there are any). Thus, if in cell $(i,j-1)$ there is a particle and cell $(i+1,j-1)$ is empty, then the particle from cell $(i,j-1)$ can go to cell $(i,j)$ or to cell $(i+1,j)$. Moreover, the choice between these possibilities can depend on what is in cells $(i,j)$, and $(i+1,j)$ (see Figure~\ref{fig:shift-up}). 

We will call this type of regulation governing the movement of particles a \emph{shifted conditional traffic-right rule} and denote it by $\Shift{}\Traffic{r}_{C}$. Note that we do not specify which shift we mean, since in this subsection it will always be the shift-up rule and we formulate the definitions below under this standing assumption. However, one can easily reformulate them in the case of any other shift rule.

\begin{definition}\label{STr}
Let $C$ be a non-empty subset of $\{00,01,10,11\}$. The shifted conditional traffic-right rule under the condition $C$ acts as follows. If $\bx\in X$ satisfies $x_{i,j-1}=1$ and $x_{i+1,j-1}=0$, then $F(\bx)_{i,j}=0$ and $F(\bx)_{i+1,j}=1$ if and only if $x_{i,j}x_{i+1,j}\in C$.
\end{definition}
In other words, if there is a particle in cell $(i,j-1)$ and cell $(i+1,j-1)$ is empty, then the particle moves to cell $(i+1,j)$ if and only if the pattern of cells $(i,j)$ and $(i+1,j)$ belongs to $C$.  

The definition of a shifted conditional traffic-left rule $\Shift{}\Traffic{l}_{C}$ is analogous.
\begin{definition}\label{STleft}
Let $C$ be a non-empty subset of $\{00,01,10,11\}$. The shifted conditional traffic-left rule under the condition $C$ acts as follows.  If $\bx\in X$ satisfies $x_{i-1,j-1}=0$ and $x_{i,j-1}=1$, then $F(\bx)_{i-1,j}=1$ and $F(\bx)_{i,j}=0$ if and only if $x_{i-1,j}x_{i,j}\in C$. \end{definition}

As outlined above, we have 15 shifted conditional traffic-right rules and 15 shifted conditional traffic-left rules. In addition, analogous to the usual conditional traffic rules, we can ``mix'' them together, obtaining four additional ones.  

In summary, if the overarching regulation is the shift-up rule, then, except for $\emptyset$, there are $34$ sets of locally applied regulations based on traffic rules: $15$ of the form $\{\Shift{}\Traffic{l}_{\Cset{}}\}$, $15$ of the form $\{\Shift{}\Traffic{r}_{\Cset{}}\}$, and four of the form $\{\Shift{}\Traffic{l}_{\Cset{}{l}}, \Shift{}\Traffic{r}_{\Cset{r}}\}$ (with $\Cset{l}$ and $\Cset{r}$ given in Eq.~(\ref{cases})). 

To conclude, let us note that the results of this and the previous subsection allow us to design 74 two-dimensional number-conserving binary CAs with the $2\times 3$-neighborhood:
\begin{itemize}
    \item Class $(\Omega,\emptyset)$ consisting of six number-conserving binary CAs with $\Lambda = \emptyset$: the identity rule and the shift rule in each of five possible directions,
    \item Class $(\Id,\Lambda)$ consisting of $34$ number-conserving binary CAs with $\Omega=\Id$ and $\Lambda\neq\emptyset$,
    \item Class $(\Shift{u},\Lambda)$ consisting of $34$ number-conserving binary CAs with $\Omega=\Shift{u}$ and $\Lambda\neq\emptyset$.
\end{itemize}

The dynamics of each of these classes is quite simple. For class $(\Omega,\emptyset)$, any local pattern will remain static (for identity) or move in one of the five directions (for shift rules). For class $(\Id,\Lambda)$, the particles can move within a given horizontal line following the conditional traffic rules, but no particle can move between the horizontal lines. Hence, the sum of states is preserved within each of the horizontal lines. For class $(\Shift{u},\Lambda)$, the dynamics is the same as for $(\Id,\Lambda)$, but the entire configuration shifts up in each time step.

In the next subsection, we prove that there are no other two-dimensional number-conserving binary CAs with the $2\times 3$-neighborhood. 

\subsection{Enumeration based on a brute force method}\label{enum}
There are as many as $2^{64}$ two-dimensional binary CAs with the $2\times 3$-neighborhood. To find all the number-conserving ones, we use the necessary and sufficient condition proved in~\cite{Durand2003} and cited below.
Note that for ease of notation we represent $f(x_{\cbi-\vv{v}{1}},x_{\cbi+\vv{v}{0}},x_{\cbi+\vv{v}{1}},x_{\cbi-\vv{v}{3}},x_{\cbi-\vv{v}{2}},x_{\cbi+\vv{v}{4}})$ graphically as 
$
\ff{x_{\cbi-\vv{v}{1}}}{x_{\cbi+\vv{v}{0}}}{x_{\cbi+\vv{v}{1}}}{x_{\cbi-\vv{v}{3}}}{x_{\cbi-\vv{v}{2}}}{x_{\cbi+\vv{v}{4}}}$, since it means simply $\ff{x_{(i-1,j)}}{x_{(i,j)}}{x_{(i+1,j)}}{x_{(i-1,j-1)}}{x_{(i,j-1)}}{x_{(i+1,j-1)}}$. 

\begin{theorem}\label{Durand}
Let $f:\{0,1\}^6\to \{0,1\}$ be the local rule of a two-dimensional binary CA with the $2\times 3$-neighborhood. This CA is number-conserving if and only if for all $x,y,z,t,u,w\in\{0,1\}$ it holds that
\begin{equation}\label{NC}
\begin{aligned}
    \ff{x}{y}{z}{t}{u}{w} =& \;x + \left[ \ff{0}{y}{z}{0}{u}{w} - \ff{0}{x}{y}{0}{t}{u}\right] + \left[ \ff{0}{0}{y}{0}{0}{u} - \ff{0}{0}{x}{0}{0}{t}\right]\\
   & + \left[ \ff{0}{0}{0}{t}{u}{w} - \ff{0}{0}{0}{x}{y}{z}\right] + \left[ \ff{0}{0}{0}{0}{y}{z} - \ff{0}{0}{0}{0}{x}{y}\right]\\
   & + \left[ \ff{0}{0}{0}{0}{t}{u} - \ff{0}{0}{0}{0}{u}{w}\right] + \left[ \ff{0}{0}{0}{0}{0}{y} - \ff{0}{0}{0}{0}{0}{x}\right]\\
   & + \left[ \ff{0}{0}{0}{0}{0}{t} - \ff{0}{0}{0}{0}{0}{u}\right] \, .
\end{aligned}
\end{equation}  
\end{theorem}

Using this condition, it is possible to check for each of the $2^{64}$ local rules whether or not it is number-conserving. However, this is quite computationally extensive. Therefore, we propose the following approach. Note that the right-hand side of Eq.~(\ref{NC}) contains only components of types $\ff{0}{\bullet}{\bullet}{0}{\bullet}{\bullet}$ and $\ff{0}{0}{0}{\bullet}{\bullet}{\bullet}$ (the symbol $\bullet$ refers to any value).
Using the inclusion–exclusion principle, we see that in order to be able to calculate the right-hand side of Eq.~(\ref{NC}), we need to define only $2^4+2^3-2^2=20$ values of $f$. Thus it is sufficient to check for each of $2^{20}$ possibilities for $\ff{0}{\bullet}{\bullet}{0}{\bullet}{\bullet}$ and $\ff{0}{0}{0}{\bullet}{\bullet}{\bullet}$ whether the value of the right-hand side of Eq.~(\ref{NC}) belongs to $\{0,1\}$.
This approach made it possible to examine all two-dimensional binary CAs with the $2\times 3$-neighborhood, and it turned out that there are $74$ of them, which is exactly the number we got in the previous subsection and which is also confirmed in~\cite{7424749}. Table~\ref{tab:empty},~\ref{tab:Id} and~\ref{tab:Su} contain all these number-conserving CAs divided into Class $(\Omega,\emptyset)$, Class $(\Id,\Lambda)$ and Class $(\Shift{u},\Lambda)$, respectively.

\begin{table}[!h]
\setlength{\arraycolsep}{8pt}
\renewcommand{\arraystretch}{0.8}
\centering
\begin{tabular}{  l  c  c  }
\toprule
 Rule number in hexagonal form\rule{0pt}{1.2em} &  \phantom{sep}$\Omega$\phantom{sep} & \phantom{sep}$\Lambda$\phantom{sep} \\
\midrule
FFFF0000FFFF0000 &  $\Id$ & $\emptyset$ \\
FF00FF00FF00FF00 &  $\Shift{l}$ & $\emptyset$ \\
FFFFFFFF00000000 &  $\Shift{r}$ & $\emptyset$ \\
CCCCCCCCCCCCCCCC &  $\Shift{u}$ & $\emptyset$ \\
AAAAAAAAAAAAAAAA &  $\Shift{\smallnwarrow}$ & $\emptyset$ \\
F0F0F0F0F0F0F0F0 &  $\Shift{\smallnearrow}$ & $\emptyset$ \\
\bottomrule
\end{tabular}
\vspace{2mm}
\caption{The list of all two-dimensional number-conserving binary CAs with the $2\times 3$-neighborhood belonging to Class $(\Omega,\emptyset)$. The first column contains the rule number being the hexadecimal representation of the binary LUT.}
\label{tab:empty}
\end{table}

\begin{table}[!h]
\setlength{\arraycolsep}{8pt}
\renewcommand{\arraystretch}{0.8}
\centering
\begin{tabular}{  l  c  c c }
\toprule
 Rule number in hexagonal form\rule{0pt}{1.2em} &  $\Omega$ & $\Cset{l}$
& $\Cset{r}$\\
\midrule
FFFF1100FCFC1100 &  $\Id$ & $\{00\}$ & - \\
FFFF2200F3F32200 &  $\Id$ & $\{01\}$ & - \\
FFFF4400CFCF4400 &  $\Id$ & $\{10\}$ & - \\
FFFF88003F3F8800 &  $\Id$ & $\{11\}$ & - \\
FFFF3300F0F03300 &  $\Id$ & $\{00,01\}$ & - \\
FFFF5500CCCC5500 &  $\Id$ & $\{00,10\}$ & - \\
FFFF99003C3C9900 &  $\Id$ & $\{00,11\}$ & - \\
FFFF6600C3C36600 &  $\Id$ & $\{01,10\}$ & - \\
FFFFAA003333AA00 &  $\Id$ & $\{01,11\}$ & - \\
FFFFCC000F0FCC00 &  $\Id$ & $\{10,11\}$ & - \\
FFFF7700C0C07700 &  $\Id$ & $\{00,01,10\}$ & - \\
FFFFBB003030BB00 &  $\Id$ & $\{00,01,11\}$ & - \\
FFFFDD000C0CDD00 &  $\Id$ & $\{00,10,11\}$ & - \\
FFFFEE000303EE00 &  $\Id$ & $\{01,10,11\}$ & - \\
FFFFFF000000FF00 &  $\Id$ & $\{00,01,10,11\}$ & - \\
FFEE0303FFEE0000 &  $\Id$ & - & $\{00\}$  \\
FFDD0C0CFFDD0000 &  $\Id$ & - & $\{01\}$  \\
FFBB3030FFBB0000 &  $\Id$ & - & $\{10\}$  \\
FF77C0C0FF770000 &  $\Id$ & - & $\{11\}$  \\
FFCC0F0FFFCC0000 &  $\Id$ & - & $\{00,01\}$ \\
FFAA3333FFAA0000 &  $\Id$ & - & $\{00,10\}$  \\
FF66C3C3FF660000 &  $\Id$ & - & $\{00,11\}$  \\
FF993C3CFF990000 &  $\Id$ & - & $\{01,10\}$  \\
FF55CCCCFF550000 &  $\Id$ & - & $\{01,11\}$  \\
FF33F0F0FF330000 &  $\Id$ & - & $\{10,11\}$  \\
FF883F3FFF880000 &  $\Id$ & - & $\{00,01,10\}$  \\
FF44CFCFFF440000 &  $\Id$ & - & $\{00,01,11\}$  \\
FF22F3F3FF220000 &  $\Id$ & - & $\{00,10,11\}$  \\
FF11FCFCFF110000 &  $\Id$ & - & $\{01,10,11\}$  \\
FF00FFFFFF000000 &  $\Id$ & - & $\{00,01,10,11\}$  \\
FF77D1C0FC741100 &  $\Id$ & $\{00\}$ & $\{11\}$ \\
FFDD2E0CF3D12200 &  $\Id$ & $\{01\}$ & $\{01\}$ \\
FFBB7430CF8B4400 &  $\Id$ & $\{10\}$ & $\{10\}$ \\
FFFF8B033F2E8800 &  $\Id$ & $\{11\}$ & $\{00\}$ \\
\bottomrule
\end{tabular}
\vspace{2mm}
\caption{The list of all two-dimensional number-conserving binary CAs with the $2\times 3$-neighborhood belonging to Class $(\Id,\Lambda)$, where $\Lambda = \{\Traffic{l}_{\Cset{l}}, \Traffic{r}_{\Cset{r}}\}$. The first column contains the rule number being the hexadecimal representation of the binary LUT.}
\label{tab:Id}
\end{table}

\begin{table}[!h]
\setlength{\arraycolsep}{8pt}
\renewcommand{\arraystretch}{0.8}
\centering
\begin{tabular}{  l  c  c c }
\toprule
 Rule number in hexagonal form\rule{0pt}{1.2em} &  $\Omega$ & $\Cset{l}$
& $\Cset{r}$\\
\midrule
CCCCCCEECCCCC0E2 &  $\Shift{u}$ & $\{00\}$ & - \\
CCCCEECCC0C0EECC &  $\Shift{u}$ & $\{01\}$ & - \\
CCEEC0C0CCEECCCC &  $\Shift{u}$ & $\{10\}$ & - \\
E2C0CCCCEECCCCCC &  $\Shift{u}$ & $\{11\}$ & - \\
CCCCEEEEC0C0E2E2 &  $\Shift{u}$ & $\{00,01\}$ & - \\
CCEEC0E2CCEEC0E2 &  $\Shift{u}$ & $\{00,10\}$ & - \\
E2C0CCEEEECCC0E2 &  $\Shift{u}$ & $\{00,11\}$ & - \\
CCEEE2C0C0E2EECC &  $\Shift{u}$ & $\{01,10\}$ & - \\
E2C0EECCE2C0EECC &  $\Shift{u}$ & $\{01,11\}$ & - \\
E2E2C0C0EEEECCCC &  $\Shift{u}$ & $\{10,11\}$ & - \\
CCEEE2E2C0E2E2E2 &  $\Shift{u}$ & $\{00,01,10\}$ & - \\
E2C0EEEEE2C0E2E2 &  $\Shift{u}$ & $\{00,01,11\}$ & - \\
E2E2C0E2EEEEC0E2 &  $\Shift{u}$ & $\{00,10,11\}$ & - \\
E2E2E2C0E2E2EECC &  $\Shift{u}$ & $\{01,10,11\}$ & - \\
E2E2E2E2E2E2E2E2 &  $\Shift{u}$ & $\{00,01,10,11\}$ & - \\
CCCCCC88CCCCFCB8 &  $\Shift{u}$ & - & $\{00\}$  \\
CCCC88CCFCFC88CC &  $\Shift{u}$ & - & $\{01\}$  \\
CC88FCFCCC88CCCC &  $\Shift{u}$ & - & $\{10\}$  \\
B8FCCCCC88CCCCCC &  $\Shift{u}$ & - & $\{11\}$  \\
CCCC8888FCFCB8B8 &  $\Shift{u}$ & - & $\{00,01\}$ \\
CC88FCB8CC88FCB8 &  $\Shift{u}$ & - & $\{00,10\}$  \\
B8FCCC8888CCFCB8 &  $\Shift{u}$ & - & $\{00,11\}$  \\
CC88B8FCFCB888CC &  $\Shift{u}$ & - & $\{01,10\}$  \\
B8FC88CCB8FC88CC &  $\Shift{u}$ & - & $\{01,11\}$  \\
B8B8FCFC8888CCCC &  $\Shift{u}$ & - & $\{10,11\}$  \\
CC88B8B8FCB8B8B8 &  $\Shift{u}$ & - & $\{00,01,10\}$  \\
B8FC8888B8FCB8B8 &  $\Shift{u}$ & - & $\{00,01,11\}$  \\
B8B8FCB88888FCB8 &  $\Shift{u}$ & - & $\{00,10,11\}$  \\
B8B8B8FCB8B888CC &  $\Shift{u}$ & - & $\{01,10,11\}$  \\
B8B8B8B8B8B8B8B8 &  $\Shift{u}$ & - & $\{00,01,10,11\}$  \\
B8FCCCEE88CCC0E2 &  $\Shift{u}$ & $\{00\}$ & $\{11\}$ \\
CCCCAACCF0F0AACC &  $\Shift{u}$ & $\{01\}$ & $\{01\}$ \\
CCAAF0F0CCAACCCC &  $\Shift{u}$ & $\{10\}$ & $\{10\}$ \\
E2C0CC88EECCFCB8 &  $\Shift{u}$ & $\{11\}$ & $\{00\}$ \\
\bottomrule
\end{tabular}
\vspace{2mm}
\caption{The list of all two-dimensional number-conserving binary CAs with the $2\times 3$-neighborhood belonging to Class $(\Shift{u},\Lambda)$, where $\Lambda = \{\Shift{}\Traffic{l}_{\Cset{l}}, \Shift{}\Traffic{r}_{\Cset{r}}\}$. The first column contains the rule number being the hexadecimal representation of the binary LUT.}
\label{tab:Su}
\end{table}

\section{The case of the two-dimensional Moore neighborhood}

The method used in the previous section easily translates to the case of any $d$-dimensional neighborhood with radius one and in this section we present how to use it for the two-dimensional Moore neighborhood. 
First of all, it can be used to design non-trivial number-conserving binary CAs with such neighborhood. 

Of course, we may start with the fact that the Moore neighborhood contains the $2\times 3$ neighborhood, so the 74 rules described in the previous section can be directly transplanted and will result in the same behavior. Additionally, since the Moore neighborhood is invariant with respect to rotation, one can define similar sets of rules that combine identity with conditional traffic-up or traffic-down rules, as well as combinations of shift-right, shift-left, and shift-down rules as overarching rule with corresponding perpendicular traffic rules. However, the additional row of cells in the neighborhood allows hundreds of thousands of new number-conserving CAs to be obtained. 

\begin{figure}[!ht]
\centering
 \subfloat[]{
\includegraphics[height=0.18\textwidth]{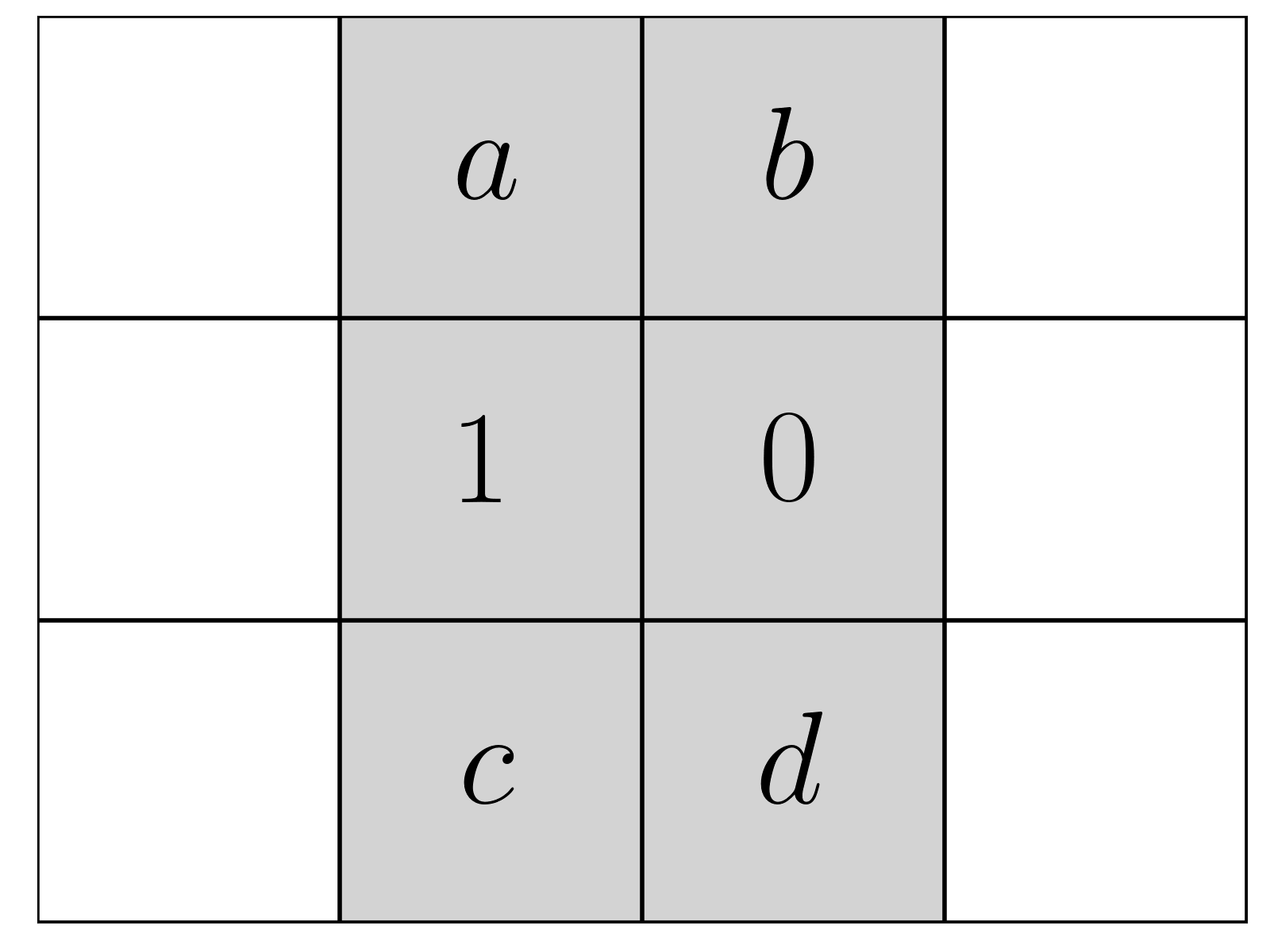}
\label{fig:TR-Moore-a}
}
 \subfloat[]{
\includegraphics[height=0.24\textwidth]{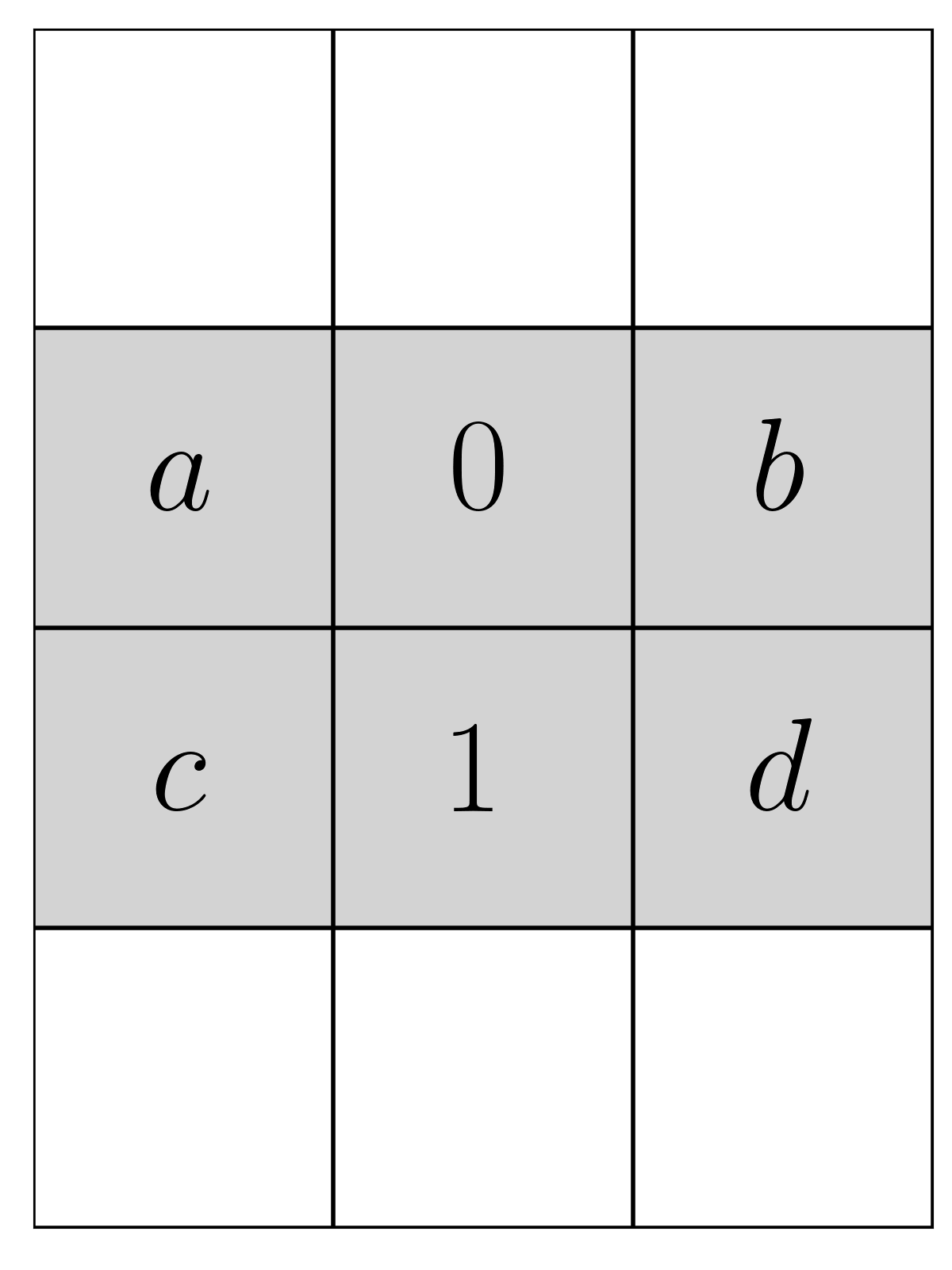}
\label{fig:TR-Moore-b}
}
 \subfloat[]{
\includegraphics[height=0.24\textwidth]{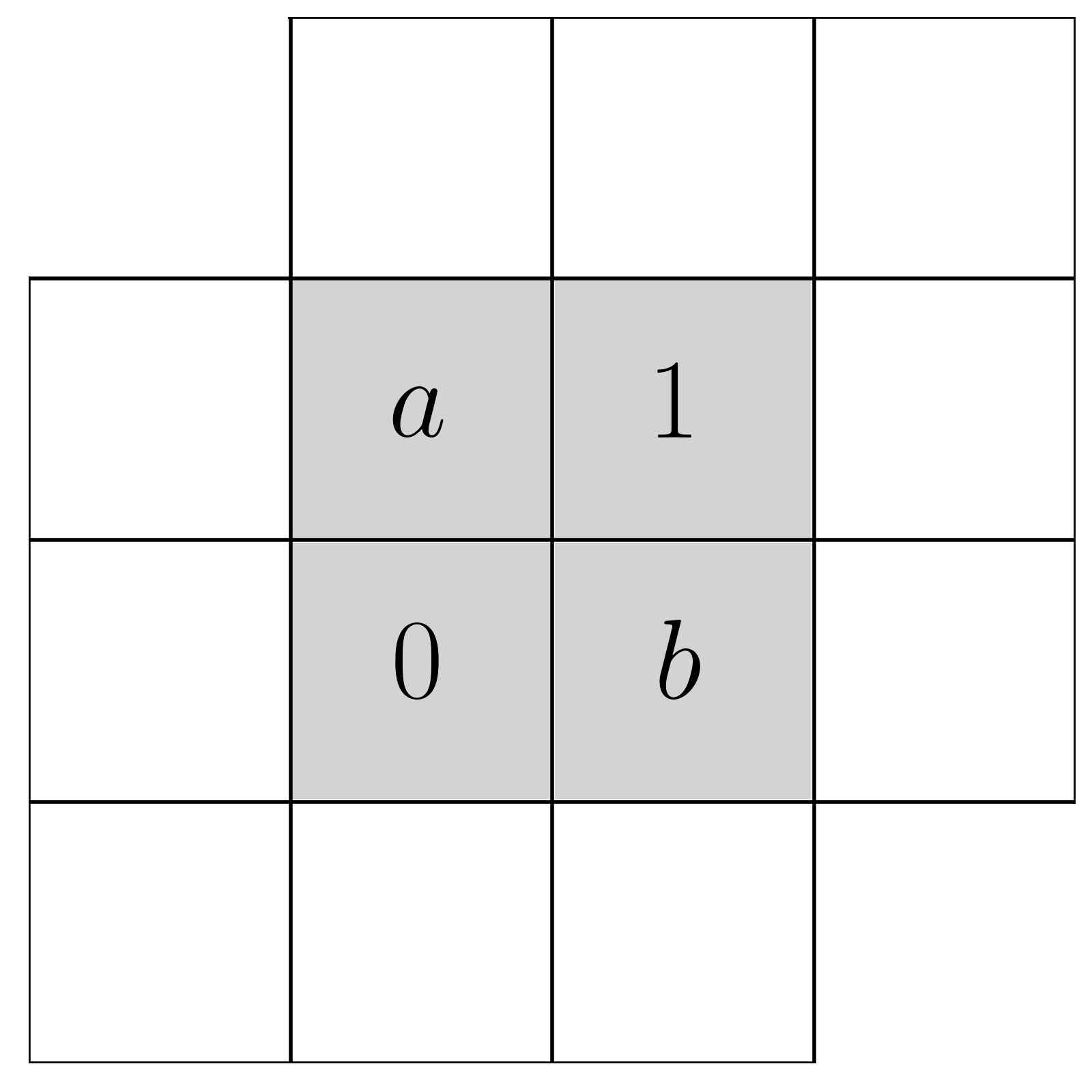}
\label{fig:TR-Moore-c}
}
\caption{Illustration of the shared neighbors in the Moore neighborhood that are relevant for the (a) traffic-right, (b) traffic-up and (c) traffic-down-left rule. Cells marked with letters $a$, $b$, $c$, and~$d$ can be used to define the corresponding conditions $\Cset{r}$, $\Cset{u}$, $\Cset{\smallswarrow}$.}
\label{fig:TR-Moore}
\end{figure}

Firstly, the conditional traffic rules may be much more complicated. For example, in the case of the conditional traffic-right rule, there are four cells visible to both the cell with a particle and its right neighbor, the patterns of which can be used as the condition (see Figure~\ref{fig:TR-Moore-a}). For this reason, there are as many as $2^{16}-1=65\,535$ conditional traffic-right rules (instead of $15$). The same is the case for the conditional traffic-left, traffic-up, and traffic-down rules.

Secondly, in addition to vertical and horizontal traffic rules, the Moore neighborhood allows us to define traffic rules along diagonal directions. These are conditional traffic-up-right, traffic-up-left, traffic-down-right, and traffic-down-left rules. For these rules, there are only two cells that can be used to define the conditions that allow traffic to occur (see Figure~\ref{fig:TR-Moore-c}), therefore there can be $15$ conditional traffic rules in each of the diagonal directions.

The last difference between the $2\times 3$-neighborhood and the Moore neighborhood is that in the latter much more complex ``mixing'' of traffic rules is possible that can lead to very complex behavior. For example, we set $\Omega = \Id$ and combine the conditional traffic rules presented in Figure~\ref{fig:TR-Moore} with the following conditions: 
\begin{itemize}
    \item The conditional traffic-right rule with condition $\begin{array}{cc}
\bullet & 0  \\
1 & \bullet
\end{array}$,\\\ie, $\Cset{r} = \{0010,0011,1010,1011\}$, illustrated by $abcd$ in Figure~\ref{fig:TR-Moore-a}.
 \item The conditional traffic-up rule with condition $\begin{array}{cc}
\bullet & 1  \\
0 & \bullet
\end{array}$,\\\ie, $\Cset{u} = \{0100,0101,1100,1101\}$, illustrated by $abcd$ in Figure~\ref{fig:TR-Moore-b}.
 \item The conditional traffic-down-left rule with condition $\begin{array}{cc}
1 &   \\
 & 0
\end{array}$,\\\ie, $\Cset{\smallswarrow} = \{10\}$, illustrated by $ab$ in 
Figure~\ref{fig:TR-Moore-c}.
\end{itemize}

It is fairly easy to verify that with this selection of conditions, requirements (R1) and (R2) are met, which ensures that the combination of these traffic rules generates a number-conserving binary CA, which, for convenience, we will refer to as CA \DD\ in the remainder of this paper.

We may now generate the LUT of the local rule of \DD, however, as it has 512 values, it is more convenient to present it in the hexadecimal representation:

\noindent
{\scriptsize FCFC0CFFFFFF0CFF FC300C00FF330C00 
 FCFCFCF000000C00 FC30FCF0FF330C00 \\
 FCFC0CFFFFFF0CFF FC300C00FF330C00 
 FCFCFCF000000C00 FC30FCF0FF330C00  } 

Let us emphasize that if we have only the LUT representation of a two-dimensional CA, our ability to study this CA would be quite limited. Of course, we can check whether it is indeed number-conserving or not, but most of the other interesting properties are undecidable. However, we can simulate the action of such a CA using a computer and formulate conclusions about its dynamics based on the obtained space-time diagrams.
\begin{figure}[!ht]
\centering
\subfloat[]{
\includegraphics[height=0.35\textwidth]{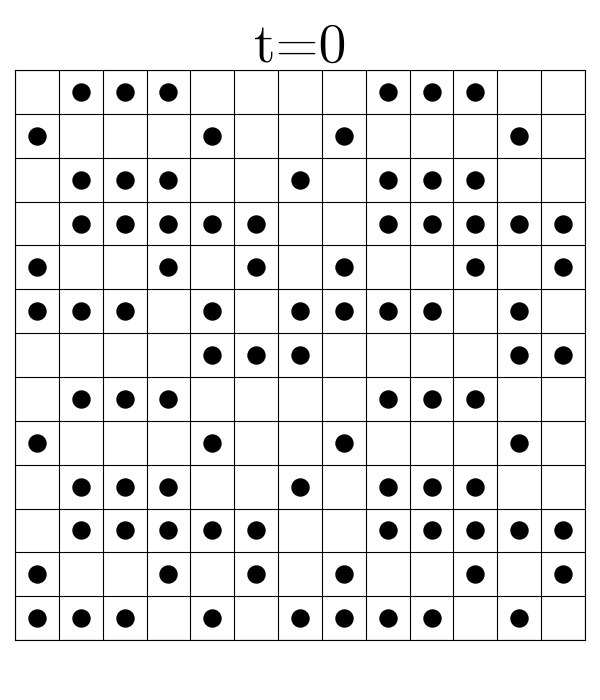}
}\quad\quad
 \subfloat[]{
\includegraphics[height=0.35\textwidth]{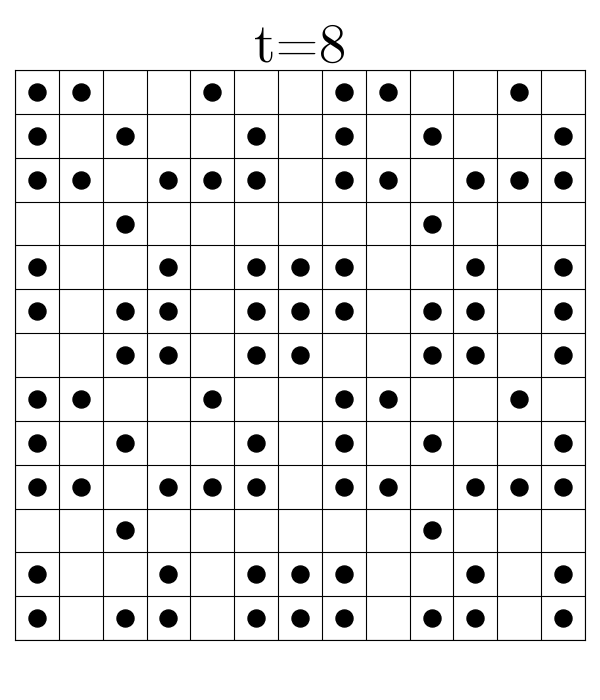}
}
\caption{A simulation of CA \DD: a fragment of a sample configuration (a) at the beginning of the simulation, (b) after 8 time steps.}
\label{fig:moore-simulation}
\end{figure}

For example, Figure~\ref{fig:moore-simulation} presents the action of CA \DD: a fragment of a sample initial configuration and this configuration after eight time steps. We observe that the behavior of \DD\ is not trivial and that its dynamics is much more complicated than that of the two-dimensional number-conserving binary CAs with the von Neumann neighborhood.

\begin{figure}[!h  ]
\centering
\subfloat[]{
\includegraphics[height=0.35\textwidth]{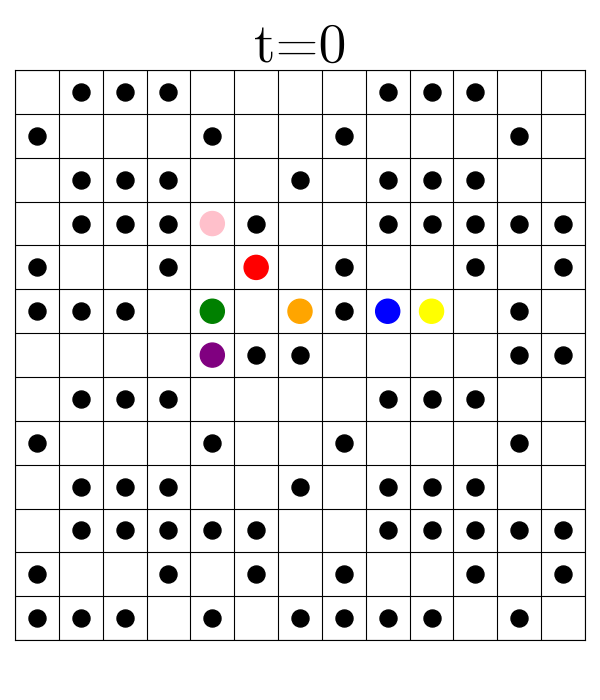}
\label{fig:moore-simulation-movement-a}
}
 \subfloat[]{
\includegraphics[height=0.35\textwidth]{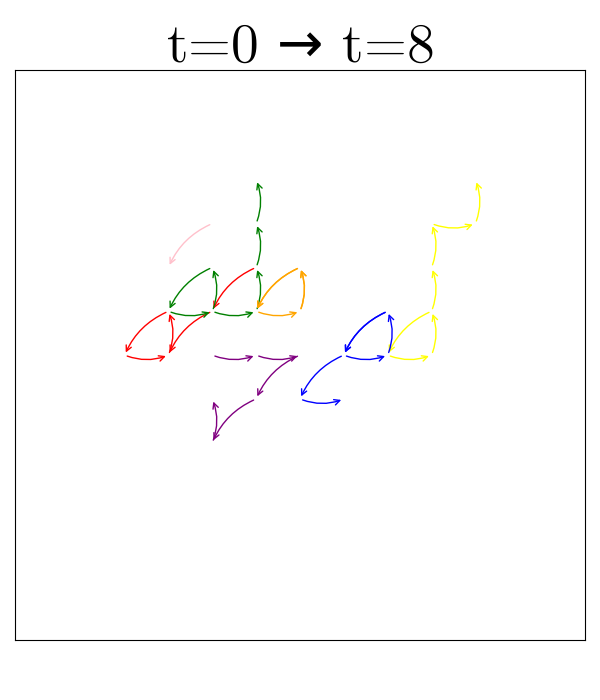}
\label{fig:moore-simulation-movement-b}
}
 \subfloat[]{
\includegraphics[height=0.35\textwidth]{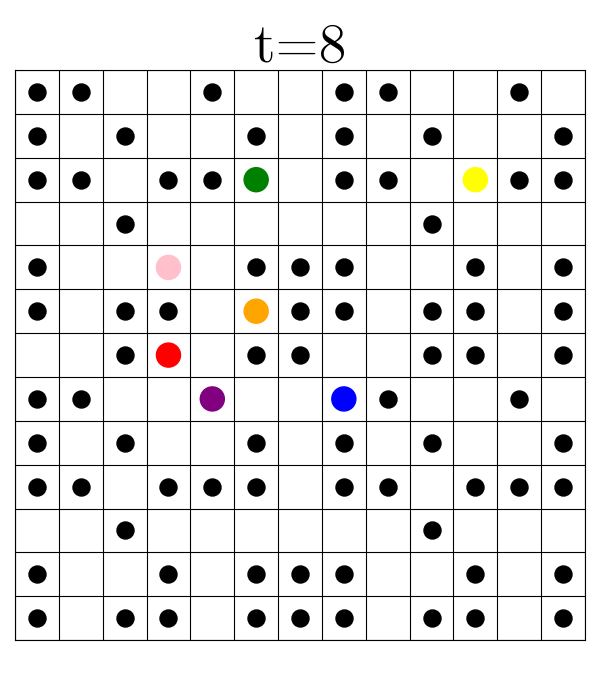}
\label{fig:moore-simulation-movement-c}
}
\caption{A simulation of CA \DD\ with the same initial configuration as shown in Figure~\ref{fig:moore-simulation}, with the trajectories of seven randomly selected particles marked. (a) The initial configuration at $t=0$, with random particles marked with colors; (b) The trajectories of the selected particles throughout the simulation; (c) The final configuration at $t=8$, with the final locations of the initially selected particles visible.}
\label{fig:moore-simulation-movement}
\end{figure}

Let us now see how representing \DD\ in terms of $\Omega$ and $\Lambda$ allows for a better understanding of its dynamics.
Indeed, since we can interpret the local rule of~\DD\ as combinations of the identity with a set of three conditional traffic rules, it is inherently easy to track the movement of each of the particles separately, as shown in Figure~\ref{fig:moore-simulation-movement}. 

Note that the first picture, Figure~\ref{fig:moore-simulation-movement-a}, presents exactly the same fragment of the initial configuration as Figure~\ref{fig:moore-simulation}, but we marked randomly selected particles with colors. Figure~\ref{fig:moore-simulation-movement-b} shows the trajectories of these particles throughout the eight steps of the simulation. Figure~\ref{fig:moore-simulation-movement-c} presents the final configuration (so it is the same as in Figure~\ref{fig:moore-simulation}), but now we can indicate the final location of each of the particles colored at the beginning.

\begin{figure}[!h  ]
\centering
\subfloat[]{
\includegraphics[height=0.30\textwidth]{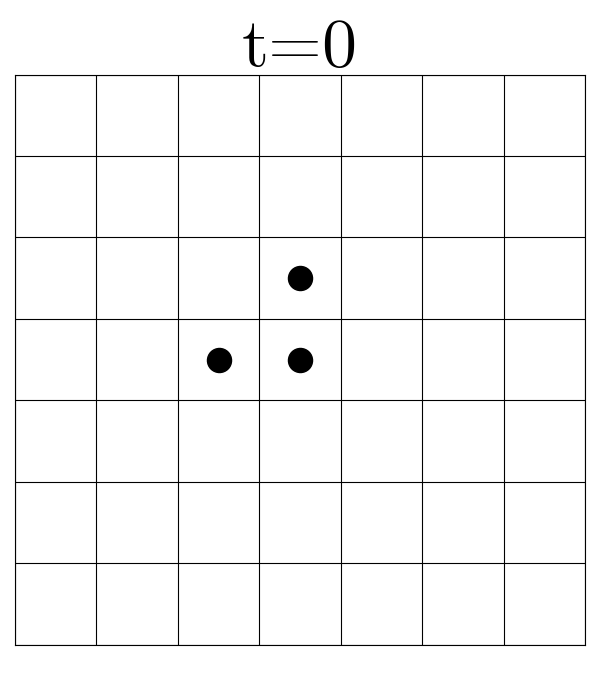}
\label{fig:example-a}
}
 \subfloat[]{
\includegraphics[height=0.30\textwidth]{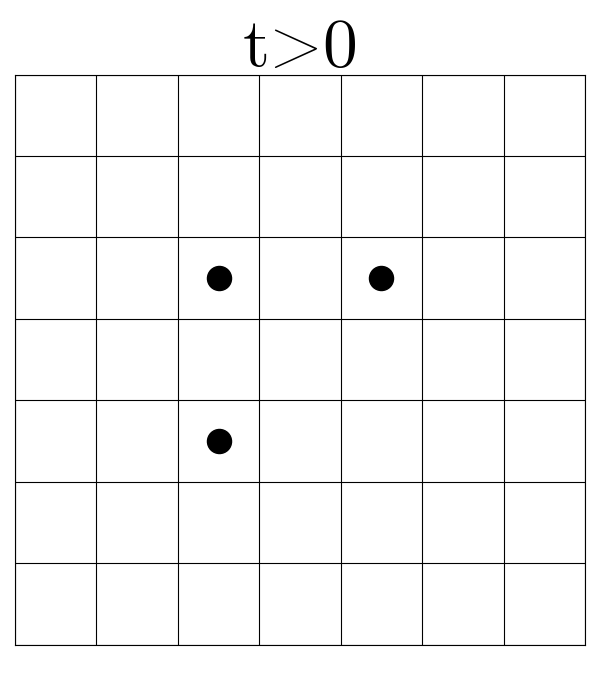}
\label{fig:example-b}
}
\caption{Illustration of the counterexample, showing that CA~\DD\ is not reversible. Starting from the initial configuration shown in (a) after one step, we reach the configuration shown in (b) which is a fixed point of CA~\DD.}
\label{fig:example}
\end{figure}

This new perspective allows us also to find the answer to whether CA~\DD\ is reversible or not. Indeed, understanding how the conditional traffic rules work, it is quite easy to design a configuration that is not a fixed point of \DD\ but is an eventually fixed point (see Figure~\ref{fig:example}), and this means that \DD\ is not injective, hence not reversible. 

\section{Conclusions and future work}
In this paper, we presented an innovative approach to describe the structure of two-dimensional binary CAs with radius one that can simulate particle flow.
By offering a simple and clear physical interpretation in terms of particle movement rules reminiscent of traffic rules, a completely new perspective on these dynamic systems is obtained. Conditional traffic rules, used in this approach, can be intuitively seen as particle flows that depend on local environmental conditions---similar to a car driver deciding to move forward not only if the lane ahead is clear, but also depending on conditions in adjacent lanes. Thus, they are simply the theoretical counterparts of realistic scenarios such as lane-changing decisions or obstacle avoidance in two-dimensional particle flows.

The presented approach, as shown in the example, allows us to design hundreds of thousands of interesting number-conserving binary CAs with radius one that have non-trivial dynamics (not only two-dimensional ones). 
Furthermore, it allows for tracking the trajectory of each individual particle throughout the entire evolution of the system. A simple description in terms of $\Omega$ and $\Lambda$ facilitates understanding of the system dynamics, which gives hope for easier verification of their key properties.

Admittedly, at this point we do not know whether CAs described in terms of a couple $(\Omega,\Lambda)$ exhaust all two-dimensional number-conserving binary CAs with the Moore neighborhood.
Therefore, the next steps in the investigation should cover at least two issues.
Firstly, it will become possible to calculate the exact number of all such couples (this rather complicated combinatorial task goes beyond the scope of this paper).
Secondly, an attempt can be made to enumerate all two-dimensional number-conserving binary CAs with the Moore neighborhood using some necessary and sufficient condition.
Clearly, the use of a condition analogous to that in Theorem~\ref{Durand} does not seem to be helpful (even when used cleverly, it still results in a  complexity of the order of $2^{112}$), but perhaps some significant improvement to the method used in Section~\ref{enum} could be helpful, or tools analogous to those created for the von Neumann neighborhood (see~\cite{decomposition}).

The introduced classification of number-conserving CAs not only constitutes theoretical advancements, but also enriches our physical understanding of particle flows simulated by CAs. Identifying CA rules as compositions of the identity or shift rules with conditional traffic rules allows these abstract models to be directly mapped onto realistic physical scenarios, and future research could further explore applications of this classification framework to model complex phenomena in fluid dynamics, granular materials, or multi-agent transport systems.

\printbibliography 

\end{document}